\documentclass{article}
\usepackage{amssymb}
\usepackage{amsmath}

\newtheorem{theorem}{Theorem}

\newtheorem{corollary}[theorem]{Corollary}

\newtheorem{definition}[theorem]{Definition}
\newtheorem{example}[theorem]{Example}

\newtheorem{lemma}[theorem]{Lemma}

\newtheorem{proposition}[theorem]{Proposition}
\newtheorem{remark}[theorem]{Remark}

\input{tcilatex}

\begin{document}

\title{On the effects of a common-pool resource on cooperation among firms
with linear technologies}
\author{\ E. Guti\'{e}rrez, N. Llorca, J. S\'{a}nchez-Soriano\thanks{%
Corresponding author, e-mail: joaquin@umh.es}\vspace*{-1mm} \\
{\small CIO and Department of Statistics, Mathematics and Computer Science}%
\vspace*{-1mm}\\
{\small University Miguel Hern\'andez of Elche, Spain}\\
\ M. A. Mosquera\\
{\small Department of Statistics and Operations Research }\\
{\small University of Vigo, Spain}}
\date{}
\maketitle

\begin{abstract}
In this paper we study the effect that the external management of a limited
(natural) resource such as carbon dioxide or water quotas has on the
behaviour of firms in a given sector. To do this, we choose a model in which
all firms have the same technology and this is lineal. In the analysis of
the problem games in partition function form arise in a natural way. It is
proved, under certain conditions, that stable allocations exist in both
cases with certainty and uncertainty.

\bigskip


\smallskip \newline

\end{abstract}

\section{Introduction}

In this paper we introduce linear production situations in which there is a
limited common-pool resource. It is managed by an external agent and is
absolutely necessary to produce any product. This type of situation appears
frequently in real-life situations related to natural resource management
such as when the producers need to buy carbon dioxide, water or fish quotas
or even to obtain public capital to invest in their firms. Imagine the case
that the producers, due to a new environmental regulation, have to introduce
a restriction concerning the emissions of greenhouse gases, for instance.
Let us assume that the greenhouse gases quotas can be bought from an
external agent, who has a limited amount, at a given market price. We wonder
what the effect of this common-pool resource on cooperation among firms will
be when they have linear production techniques.

Linear production $(LP)$ situations are situations where several producers
own resource bundles. They can use these resources to produce various
products via linear production techniques that are available to all
producers. The goal of each producer is to maximize their profit, which
equals the revenue of their products at the given market prices. These
situations and corresponding cooperative games are introduced in Owen
(1975), where it is shown that these games always have a non empty core by
constructing a core-element via a related dual linear program. Gellekom et
al. (2000) named the set of all the core-elements that can be found in this
way, the `Owen set', and they provide a characterization of this. More
general are situations involving a countably infinite number of products
that can be produced. Tijs et al. (2001) study relations between the Owen
set and the core of these semi-infinite situations.

When a limited common-pool resource is introduced in an $LP$ situation this
leads to a linear production situation with a common-pool resource $(LPP$
situation$).$ Although it may intuitively seem that\ when introducing a
small change in the $LP$ model everything will work in a similar way, in
this case it is not true because, for instance, the games that arise in
these situations are partition function form games and the existence of
stable allocations is not always guaranteed.

Models with a limited external common resource have generally been
addressed, in the literature, from the non-cooperative perspective. However,
as Hardin (1968) points out the so-called tragedy of commons, where the
common pool resource is overused, can occur. This is the main reason why we
have considered a cooperative perspective in our model. Driessen and
Meinhardt (2001) use a classic cooperative game, defined from a
non-cooperative one, in which the value of a coalition (group of producers)
is obtained from a two-person game where the members of the coalition try to
maximize their profits in the worst case. Funaki and Yamato (1999) provide a
cooperative approach to a model of an economy with a common-pool resource,
where the demands on it are additive. In our case, we are interested in
addressing $LP$ situations with a limited common-pool resource from a
cooperative point of view, where the value of a coalition is determined by
taking into account not only what the members of the coalition can do, but
also what outsiders can do. Thus, our model implies the use of games in
partition function form introduced by Thrall and Lucas (1963) as in Funaki
and Yamato (1999). In their model the distribution of fish (common-pool
resource) among fishermen is carried out in proportion to the amount of
labor involved.

This paper tackles the problem of allocating the common-pool resource
focused on $LPP$ situations. As in Funaki and Yamato (1999), the core of the
games in partition function form associated with these situations can be
reduced to the core of a related game in characteristic function form. The
amount of the common-pool resource available can play a crucial role in the
analysis of the associated games. We distinguish two cases: when the
cooperation of all producers enables the common-pool resource to be
sufficient for them and when it is not sufficient. In the former case, we
show that the core of the games in partition function form is non empty. In
the latter situation, additional conditions are needed to assure the non
emptiness of the core. The analysis of both cases is carried out taking into
account that the partition function is unknown because the producers do not
know how the common-pool resource is to be assigned. Therefore, this problem
is approached as one with uncertainty. The study, where the partition
function is known exactly, is left for further reseach.

The paper is organized as follows. Section 2 contains basic concepts on
cooperative transferable utility games. In Section 3 linear production
situations with a common-pool resource are introduced. We show that if the
common-pool resource is not a constraint for the production process, or if
it is so only for the coalition of all producers, these games can be reduced
to games in characteristic function form. In the first case, we can find
allocations in the core of the corresponding games, but in the second one
the core can be empty. We introduce a new concept for partitions to be
partitionally stable, which allows us, in some sense, to extend the concept
of the core. In Section 4 we assume that the producers do not know how the
common-pool resource will be distributed. Therefore, we introduce
common-pool resource games to deal with this uncertainty. Different points
of view can be used to define these games, we study the two extremes: the
optimistic and pessimistic. Section 5 concludes.

\section{Preliminaries}

Let $N$ be a non empty finite set of agents who agree to coordinate their
actions. A cooperative game in characteristic function form is an ordered
pair $\left( N,v\right) ,$ where $N$ is the set of players and $%
v:2^{N}\rightarrow \mathbb{R}$ is the characteristic function with $%
v(\emptyset )=0$. This function assigns to each group of players
(coalition), $S\subset N,$ the value $v(S)$ which represents what the
members in $S$ obtain when they cooperate jointly. In a transferable utility
game ($TU$-game) it is assumed that the utility can be linearly transferred
among agents.

A classic issue in cooperative game theory is how to distribute the profit
generated by the cooperating players. One way to do this is to use
allocations in the core of the game. The core, $C\left( v\right) ,$ of a
characteristic function form game $\left( N,v\right) ,$ introduced by
Gillies (1953), is the subset of vectors in $\mathbb{R}^{N}$ satisfying 
\begin{equation*}
\begin{array}{ll}
\text{(Efficiency)} & \sum_{i\in N}x_{i}=v(N)\text{, and} \\ 
\text{(Coalitional rationality)} & \sum_{i\in S}x_{i}\geq v(S),\text{ for
all }S\subset N\text{.}%
\end{array}%
\end{equation*}

Bankruptcy problems were first introduced by O'Neill (1982). A standard
bankruptcy problem can be described by a triple $(N,E,d)$, where $N=\left\{
1,...,n\right\} $ is the finite set of agents, $E\geq 0$ is the estate to be
divided and $d\in \mathbb{R}_{+}^{N},$ the vector of claims, is such that $%
\sum_{i\in N}d_{i}\geq E$. To deal with a bankruptcy problem $(N,E,d)$ we
can derive a classical bankruptcy $TU$-game $(N,v)$, where the value of a
coalition $S\subset N$ is given by 
\begin{equation*}
v(S)=\max \{E-\sum_{i\in N\setminus S}d_{i},0\},
\end{equation*}%
and represents what is left for the players in $S$ after the demands of the
players in $N\backslash S$ have been satisfied. These games have a non empty
core.

In order to deal with our model, we will need to consider games in partition
function form introduced by Thrall and Lucas (1963), where the worth of a
coalition $S$\ depends not only on what the members in $S$ can do, but also
on what outsiders do. These are cooperative games with externalities.
Formally, let $\mathcal{P}(N)$ denote the set of all partitions of $N$ and $%
P=\{S_{1},\ldots ,S_{k}\}$ represents one of these partitions or coalition
structures, where the coalitions $S_{1},\ldots ,S_{k}$\ are disjoint and
their union is $N$. The pair $\left( \left. S\right\vert P\right) $ such
that $S\in P$\ is usually called an embedded coalition. A cooperative game
in partition function form is defined by $\left( N,\mathcal{P}(N),\left\{
V\left( \left. \bullet \right\vert P\right) \right\} _{P\in \mathcal{P}%
(N)}\right) $, where $N$ is the set of players, $\mathcal{P}(N)$ denotes the
set of all partitions of $N$ and $V\left( \left. S\right\vert P\right) $
with $S\in P$\ is a real number that represents the profit that a coalition $%
S\subseteq N$ can obtain when $P$ is formed. Note that the profit that a
coalition can obtain depends on the coalitions formed by the other players
in $P\in \mathcal{P}(N)$.

Given a partition $P\in \mathcal{P}(N)$, a vector $x\in \mathbb{R}^{n}$ is
said to be feasible under $P$ if it satisfies $\tsum\limits_{i\in
S}x_{i}\leq V\left( \left. S\right\vert P\right) ,\forall S\in P.$ We denote
by $\mathcal{F}^{P}$ the set of all feasible vectors under $P$ and $\mathcal{%
F}=\cup _{P\in \mathcal{P}(N)}\mathcal{F}^{P}$ denotes the set of all
feasible vectors. Given two vectors $x,x^{\prime }$ in $\mathbb{R}^{n}$, as
in Funaki and Yamato (1999) we say that $x$ dominates $x^{\prime }$ through $%
S$ and denote $x~dom_{S}~x^{\prime }$ if the following conditions are
satified:

\begin{enumerate}
\item $\tsum\limits_{i\in S}x_{i}\leq V\left( \left. S\right\vert P\right)
,\forall P\in \mathcal{P}(N)$ such that $S\in P,$

\item $x_{i}>x_{i}^{\prime },\forall i\in S.$
\end{enumerate}

We say that $x$ dominates $x^{\prime }$ if there exists $S\subseteq N$ such
that $x~dom_{S}~x^{\prime },$ and denote $x~dom~x^{\prime }$. The core of a
cooperative game in partition function form is defined by $C\left( V\right)
=\left\{ x\in \mathcal{F}\left\vert \nexists x^{\prime }\in \mathcal{F}\text{
s.t. }x^{\prime }~dom~x\right. \right\} .$ However, if we consider another
definition of dominance, then we will obtain a different core. Thus, if we
change condition 1 by 
\begin{equation*}
\overline{1}.\tsum\limits_{i\in S}x_{i}\leq V\left( \left. S\right\vert
P\right) ,\text{ for some }P\in \mathcal{P}(N)\text{ with }S\in P,
\end{equation*}%
we obtain a more restrictive concept of dominance that we denote by $%
\overline{dom}$ and the corresponding core is defined as $\overline{C}\left(
V\right) =\left\{ x\in \mathcal{F}\left\vert \nexists x^{\prime }\in 
\mathcal{F}\text{ s.t. }x^{\prime }~\overline{dom}~x\right. \right\} .$

Associated with each game in partition function form two cooperative games
in characteristic function form can be introduced: $\left( N,v^{-}\right) $
and $\left( N,v^{+}\right) $, where%
\begin{equation*}
\begin{array}{l}
v^{-}\left( S\right) =\min \left\{ V\left( \left. S\right\vert P\right)
\left\vert P\in \mathcal{P}(N)\text{ such that }S\in P\right. \right\} , \\ 
v^{+}\left( S\right) =\max \left\{ V\left( \left. S\right\vert P\right)
\left\vert P\in \mathcal{P}(N)\text{ such that }S\in P\right. \right\} .%
\end{array}%
\end{equation*}

$\left( N,v^{-}\right) $ represents a pessimistic point of view of the gain
that a coalition $S$ can get, while $\left( N,v^{+}\right) $\ can be seen as
optimistic. Funaki and Yamato (1999) proved that if $V\left( \left. \left\{
N\right\} \right\vert N\right) >\tsum\limits_{S\in P}V\left( \left.
S\right\vert P\right) ,\forall P\in \mathcal{P}(N),$ then

\begin{description}
\item[$a)$] $C(V)=C(v^{-}),$ and

\item[$b)$] $\overline{C}\left( V\right) =C(v^{+}).$
\end{description}

Given $P,P^{\prime }\in \mathcal{P}(N)$, $P^{\prime }$ is a refinement of $P$
if for all $S^{\prime }\in P^{\prime }$ there exits $S\in P$ such that $%
S^{\prime }\subseteq S$, and it is denoted by $P^{\prime }\subseteq P$.
Using the concept of refinement an ordering of partitions arise in a natural
way, with this ordering $\left( \mathcal{P}(N),\subseteq \right) $ is the
so-called partition lattice.

\section{Linear production situations with a common-pool resource}

Let $N=\{1,\ldots ,n\}$ be a set of producers that face a linear production
problem to produce a set $G=\{1,\ldots ,g\}$ of goods from a set $%
Q=\{1,\ldots ,q\}$ of resources. There exists an external common-pool
resource, limited by an amount of $r,$ that agents need to buy for producing
the goods. The parameters of the model are:

\begin{itemize}
\item $b^{i}\in \mathbb{R}_{+}^{q}$ are the available resources for each
producer $i\in N$, $b^{S}=\sum_{i\in S}b^{i}$. $B\in \mathcal{M}_{q\times n}$
is the resource matrix. We assume that there is a positive quantity
available of each resource, that is, for all resources $t\in Q$ there is a
producer $i$ such that $b_{i}^{t}>0$.

\item The common pool-resource is not endowed to the producers but managed
by an external agent. Its cost per unit is $c$ and the total available is
denoted by $r$.

\item $A\in \mathcal{M}_{\left( q+1\right) \times g}$ is the production
matrix, $a_{tj}$ represents the amount of the resource $t$ needed to produce
the product $j$, where the last row is related to the common-pool resource
and $a_{(q+1)j}>0\quad \forall \,j\in G$. Furthermore, we do not allow for
output without input\ and therefore there exists at least one resource $t\in
Q$ with $a_{tj}>0\quad \forall \,j\in G$.

\item $p\in \mathbb{R}_{++}^{g}$ is the price vector. Moreover, in order to
deal with a profitable process we assume that $p_{j}>a_{(q+1)j}c\quad
\forall \,j\in G$.
\end{itemize}

Therefore, a linear production situation with a common-pool resource ($LPP$%
)\ can be represented by $\left( A,B,p,r,c\right) .$

To maximize his profit, producer $i$ needs an optimal production plan $%
(x;z)\in \mathbb{R}_{+}^{g+1}$ that tells him how much he should produce of
each good, $x,$ and how much he needs of the common-pool resource, $z$. Not
all production plans are feasible since the producer has to take into
account his limited amount of resources. The amount of resources needed in a
feasible production plan should not exceed the amount of resources
avaliables for producer $i$. Furthermore, a feasible production plan has to
be nonnegative since we are only interested in producing nonnegative
quantities of the products. The following linear program maximizes the
profit of producer $i$%
\begin{equation*}
\begin{array}{ll}
\max & \sum_{j=1}^{g}p_{j}x_{j}-cz \\ 
\text{s.t:} & Ax\leq \left( 
\begin{array}{c}
b^{i} \\ 
z%
\end{array}%
\right) \\ 
& x\geq \mathbf{0}_{g},z\geq 0.%
\end{array}%
\end{equation*}%
Thus, an optimal production plan for producer $i$ is an optimal solution of
this linear program. Apart from producing on their own, producers are
allowed to cooperate. If a coalition $S$ of producers cooperates then they
put all their resources together and so given this amount of resources, the
coalition wishes to maximize its profit,%
\begin{equation}
\begin{array}{ll}
\max & \sum_{j=1}^{g}p_{j}x_{j}-cz \\ 
\text{s.t:} & Ax\leq \left( 
\begin{array}{c}
b^{S} \\ 
z%
\end{array}%
\right) \\ 
& x\geq \mathbf{0}_{g},z\geq 0.%
\end{array}
\label{eq:LP2}
\end{equation}

With an abuse of notation, we use $z$ to represent the amount of the
common-pool resource that a producer or a group of producers will need. We
denote by $value\left( S;z\right) $ the value of this linear program, for
every fixed $z$.

The optimal demand of the common-pool resource for each coalition $S,$ $%
d_{S}=\min \left\{ z\in \mathbb{R}_{+}\left\vert value\left( S;z\right) 
\text{ is maximum}\right. \right\} ,$ is obtained by solving the linear
program (\ref{eq:LP2}). We should point out that these optimal demands are
the desired amount of the common-pool resource for each coalition $S$ and
can be seen as their utopic or greatest aspirations a priori, before the
common-pool resource is allocated. Note that they are not bounded from above
by $r.$

Although it may seem that these demands are superadditive, i. e. $d_{S}\geq
\sum_{i\in S}d_{\{i\}},$ this is not true as the next example shows.

\begin{example}
\label{ex. super}Let $(A,B,r,p,c)$ be an $LPP$ situation, with two
producers, $N=\left\{ 1,2\right\} ,$ who produce three products from two
resources and a common-pool resource, where 
\begin{equation*}
A=\left[ 
\begin{array}{lll}
1 & 0 & 1 \\ 
0 & 1 & 1 \\ 
2 & 2 & 1%
\end{array}%
\right] ,B=\left[ 
\begin{array}{ll}
4 & 1 \\ 
1 & 4%
\end{array}%
\right] ,p=\left( 
\begin{array}{l}
4 \\ 
4 \\ 
8%
\end{array}%
\right) ,c=1,r=5.
\end{equation*}%
In this case, $d_{\left\{ 1\right\} }=d_{\left\{ 2\right\} }=7$ while $%
d_{\left\{ 12\right\} }=5.$
\end{example}

Next we present a technical result which guarantees that once we know that a
positive profit is achieved, all the lower levels of the common-pool
resource also provide positive profits.

\begin{proposition}
Let $S\subseteq N$, if there is $z^{\ast }$ such that $value\left( S;z^{\ast
}\right) >0$, then $value\left( S;z\right) >0,$ for all $z$ with $%
0<z<z^{\ast }.$
\end{proposition}

\noindent \textbf{Proof. }Let $z$ such that $0<z<z^{\ast }$ and $x^{\ast }$
the optimal solution corresponding to $value\left( S;z^{\ast }\right) .$
Consider $x^{z}=\frac{z}{z^{\ast }}x^{\ast }$. The point $(x^{z};z)$ is
feasible for the problem corresponding to $S$. Furthermore, it holds that:%
\begin{equation*}
\begin{array}{ll}
\sum_{j=1}^{g}p_{j}x_{j}^{\ast }-cz^{\ast }>0\Longrightarrow & \frac{z}{%
z^{\ast }}\left( \sum_{j=1}^{g}p_{j}x_{j}^{\ast }-cz^{\ast }\right)
>0\Longrightarrow \\ 
\sum_{j=1}^{g}p_{j}\left( \frac{z}{z^{\ast }}x_{j}^{\ast }\right) -c\left( 
\frac{z}{z^{\ast }}z^{\ast }\right) >0\Longrightarrow & 
\sum_{j=1}^{g}p_{j}x_{j}^{z}-cz>0.%
\end{array}%
\end{equation*}%
\hfill \hfill\ \hfill .\hfill $\blacksquare $\bigskip

In the sequel, we will assume that for all $S$, there is a feasible
production plan $(x;z)$ such that $value\left( S;z\right) >0.$ This implies
that $d_{S}>0.$

Let us assume that $P$ is formed and, either through a collaborative
procedure or through a competitive mechanism\footnote{%
At this moment we do not specify any process, what follows holds in any case.%
}$,$ the amount of the common-pool resource\ finally allocated to coalition $%
S\in P$ by the manager is $z_{S}(P).$ The profit that a coalition $%
S\subseteq N$ can obtain is given by 
\begin{equation}
\begin{array}{ll}
\max & \sum_{j=1}^{g}p_{j}x_{j}-cz_{S}(P) \\ 
\text{s.t:} & Ax\leq \left( 
\begin{array}{c}
b^{S} \\ 
z_{S}(P)%
\end{array}%
\right) \\ 
& x\geq \mathbf{0}_{g}.%
\end{array}
\label{eqLP1}
\end{equation}

We should point out that $z_{S}(P)$ is bounded from above by $r$, because
the manager of the common-pool resource cannot exceed this amount, while
this does not hold for $d_{S},$ for all $S\subseteq N.$

Depending on the procedure used to obtain $z_{S}(P),$ which will be less or
equal to its optimal demand $d_{S}$, we can define different games. These
games are not characteristic function form games, but partition function
form games and the amount available of the common-pool resource can play a
crucial role on analysis.

\begin{definition}
Let $\left( A,B,p,r,c\right) $ be an $LPP$\ situation. The partition
function form game associated with this situation is given by $\left( N,%
\mathcal{P}(N),\left\{ V\left( \left. \bullet \right\vert P\right) \right\}
_{P\in \mathcal{P}(N)}\right) $, where $N$ is the set of players, $\mathcal{P%
}(N)$ denotes the set of all partitions of $N$ and $V\left( \left.
S\right\vert P\right) $ with $S\in P$ is obtained from the linear program (%
\ref{eqLP1}), for all $S\subset N,$ where $z_{S}(P)$ is the amount of
common-pool resource avaliable for coalition $S$ when partition $P$ is
formed.
\end{definition}

\begin{proposition}
Let $\left( A,B,p,r,c\right) $ be an $LPP$\ situation and $\left( N,\mathcal{%
P}(N),\left\{ V\left( \left. \bullet \right\vert P\right) \right\} _{P\in 
\mathcal{P}(N)}\right) $\ the corresponding partition function form game.
Then, 
\begin{equation*}
V\left( \left. N\right\vert \left\{ N\right\} \right) \geq
\tsum\limits_{S\in P}V\left( \left. S\right\vert P\right) ,\forall P\in 
\mathcal{P}(N).
\end{equation*}
\end{proposition}

\noindent \textbf{Proof. }Given $P\in \mathcal{P}(N),V\left( \left.
S\right\vert P\right) =value\left( S;z_{S}(P)\right) ,\forall S\in P$ such
that $\tsum\limits_{S\in P}z_{S}(P)\leq r.$ Let be $\left(
x^{S};z_{S}(P)\right) $ an optimal plan for each coalition $S\in P$. Thus, $%
Ax^{S}\leq \left( 
\begin{array}{c}
b^{S} \\ 
z_{S}(P)%
\end{array}%
\right) $ and 
\begin{equation*}
A\left( \tsum\limits_{S\in P}x^{S}\right) \leq \left( 
\begin{array}{c}
\tsum\limits_{S\in P}b^{S} \\ 
\tsum\limits_{S\in P}z_{S}(P)%
\end{array}%
\right) \leq \left( 
\begin{array}{c}
b^{N} \\ 
r%
\end{array}%
\right) .
\end{equation*}%
Then, $\left( \tsum\limits_{S\in P}x^{S};\tsum\limits_{S\in
P}z_{S}(P)\right) $ is a feasible production plan for $N$ and 
\begin{equation*}
\tsum\limits_{S\in P}value\left( S;z_{S}(P)\right) \leq value\left(
N;\tsum\limits_{S\in P}z_{S}(P)\right) \leq V\left( \left. N\right\vert
\left\{ N\right\} \right) .
\end{equation*}%
\hfill $\blacksquare $\bigskip

It is easy to check that $C\left( V\right) $ and $\overline{C}\left(
V\right) $ only include efficient allocations. The following consequences
are give without a proof because they can be derived in a similar manner as
in Funaki and Yamato (1999).

\begin{corollary}
Let $\left( A,B,p,r,c\right) $ be an $LPP$\ situation, $\left( N,\mathcal{P}%
(N),\left\{ V\left( \left. \bullet \right\vert P\right) \right\} _{P\in 
\mathcal{P}(N)}\right) $\ the corresponding partition function form game and 
$\left( N,v^{-}\right) ,$ $\left( N,v^{+}\right) $ the related games in
characteristic function form. Then, $C\left( V\right) =C\left( v^{-}\right) $
and $\overline{C}\left( V\right) =C\left( v^{+}\right) .$
\end{corollary}

If $\left( N,\mathcal{P}(N),\left\{ V\left( \left. \bullet \right\vert
P\right) \right\} _{P\in \mathcal{P}(N)}\right) $ is such that$,\forall
S\subseteq N,\forall P\in \mathcal{P}(N)$ with $S\in P$, $V\left( \left.
S\right\vert P\right) =v\left( S\right) $, then the two definitions of
dominance are equivalent, $\left( N,v\right) $ is a game in characteristic
function form and the core reduces to the well-known definition for games in
characteristic function form. Next we show other two situations in which
this holds.

Given a partition $P=\{S_{1},\ldots ,S_{k}\},$ its total demand is $%
d(P)=\sum_{i=1}^{k}d_{S_{i}}$. An outstanding set associated with both the
common-pool resource and the set of all partitions is the following:%
\begin{equation*}
M^{\text{min}}=\{P\in \mathcal{P}(N):\ d(P)>r\text{ and }P\text{ is minimal
for the operator }\subseteq \},
\end{equation*}
where $P$ is minimal for the operator $\subseteq $ means that there is no
refinement $P^{\prime }$ of $P$ with $d(P^{\prime })>r.$

The next results show that if the common-pool resource is not a constraint
on the production process, i.e. $M^{\text{min}}=\varnothing ,$\ or it is
only a restriction for the grand coalition, i.e. $M^{\text{min}}=\left\{
N\right\} ,$\ then the corresponding partition function form games are
characteristic function form games.

\begin{proposition}
\label{max}Let $\left( A,B,p,r,c\right) $ be an $LPP$\ situation. If $M^{%
\text{min}}=\varnothing ,$ then the corresponding game $\left( N,\mathcal{P}%
(N),\left\{ V\left( \left. \bullet \right\vert P\right) \right\} _{P\in 
\mathcal{P}(N)}\right) $ is a characteristic function form game.
\end{proposition}

\noindent \textbf{Proof. }If $M^{\text{min}}=\varnothing ,$ then we have $%
d_{N}\leq r.$ Thus, $z_{S}(P)=d_{S}$ for all $S\subseteq N,$ due to $d\left(
P\right) \leq r$ for all $P$. Therefore, for each $S\subseteq N,$ $V\left(
\left. S\right\vert P\right) =V\left( \left. S\right\vert P^{\prime }\right) 
$ for all $P,P^{^{\prime }}\in \mathcal{P}(N)$ such that $S\in P$, i.e. for
each coalition the value does not depend on the coalitions formed by other
players. \hfill $\blacksquare $\bigskip

\begin{proposition}
\label{N>r}Let $\left( A,B,p,r,c\right) $ be an $LPP$ situation. If $\left\{
N\right\} \in M^{\text{min}},$ then $M^{\text{min}}=\left\{ N\right\} $ and
the related game $\left( N,\mathcal{P}(N),\left\{ V\left( \left. \bullet
\right\vert P\right) \right\} _{P\in \mathcal{P}(N)}\right) $ is a
characteristic function form game.
\end{proposition}

\noindent \textbf{Proof. }If $\left\{ N\right\} \in M^{\text{min}}$ then,
since $\left\{ N\right\} $ is minimal for the operator $\subseteq ,$ $%
d_{N}>r $ and $d\left( P\right) \leq r$ for all $P\in $ $\mathcal{P}(N)$, $%
P\neq \left\{ N\right\} .$ Similarly to Proposition \ref{max}, for each $%
S\subseteq N,V\left( \left. S\right\vert P\right) =V\left( \left.
S\right\vert P^{\prime }\right) $ for all $P,P^{^{\prime }}\in \mathcal{P}%
(N),$ i.e. the value of coalition $S$\ does not depend on the coalitions
formed by other players. On the other hand, the value of the grand coalition 
\begin{equation}
\begin{array}{ll}
\max & \sum_{j=1}^{g}p_{j}x_{j}-cr \\ 
\text{s.t:} & Ax\leq \left( 
\begin{array}{c}
b^{N} \\ 
r%
\end{array}%
\right) \\ 
& x\geq \mathbf{0}_{g}.%
\end{array}
\label{eq:LP3}
\end{equation}%
only depends on its own, since there is no partition including $N$ as a
proper subset. Therefore, the related game is a characteristic function form
game.\hfill $\blacksquare $\bigskip

Let $\left( A,B,p,r,c\right) $ be an $LPP$ situation. The characteristic
function form game associated with one of the two previous situations $%
\left( N,v\right) $, where $d_{N}\leq r$ or $d_{N}>r$ and $d\left( P\right)
\leq r$ for all $P\neq \left\{ N\right\} ,$ is given by $v\left( S\right)
=value\left( S;z\right) $, with $z=d_{S}$ for all $S\neq N$ and $z=\min
\left\{ d_{N},r\right\} $ for the grand coalition $N.$ This is due to the
fact that the common-pool resource is sufficient to satify the demands for
all $S\neq N,$ but for $N$ the maximum amount avaliable is $r$.

The next result shows that the characteristic function form game obtained
when the common-pool resource is not a constraint for the production process
has a non empty core.

\begin{theorem}
\label{dual}Let $\left( A,B,p,r,c\right) $ be an $LPP$ situation with $M^{%
\text{min}}=\varnothing $. The characteristic function form game $\left(
N,v\right) $ associated with this situation has a non empty core.
\end{theorem}

\noindent \textbf{Proof. }The dual problem of (\ref{eq:LP2})\footnote{%
We use this problem because it is known by hypotesis that $\exists z\leq r$
for all problems, i.e. the common-pool resource is not scarce in any case.}
for the grand coalition, $N$, is 
\begin{equation}
\begin{array}{ll}
\min & \sum_{t=1}^{q}b_{t}^{N}y_{t}+0y_{q+1} \\ 
\text{s.t:} & A^{t}y\geq p \\ 
& y_{q+1}\leq c \\ 
& y\geq \mathbf{0}_{q+1}.%
\end{array}
\label{dual N}
\end{equation}%
An optimal solution of (\ref{eq:LP2}) for the grand coalition $N$ is given
by $\left( x^{N};d_{N}\right) $ with $d_{N}\leq r,$ and the related dual
optimal solution is $\left( y_{q}^{N};y_{q+1}^{N}\right) $, where with an
abuse of notation from now on, we represent by $y_{q}^{N}$ the vector $%
\left( y_{1}^{N},...,y_{q}^{N}\right) $. From duality, it is known that $%
\sum_{j=1}^{g}p_{j}x_{j}^{N}-cd_{N}=%
\sum_{t=1}^{q}b_{t}^{N}y_{t}^{N}+0y_{q+1}^{N}=v\left( N\right) .$ Therefore,
somehow, the cost of the common-pool resource is charged to (discounted
from) the value of the resources. It is easy to check that $\left(
y_{q}^{N};y_{q+1}^{N}\right) $ is feasible in the dual problem of (\ref%
{eq:LP2}) for every coalition $S\subset N.$ Moreover, we have that for a
dual optimal solution $\left( y_{q}^{S};y_{q+1}^{S}\right) $ associated with
the optimal solution $\left( x^{S};d_{S}\right) ,$ it holds that $%
\sum_{t=1}^{q}b_{t}^{S}y_{t}^{N}+0y_{q+1}^{N}\geq
\sum_{t=1}^{q}b_{t}^{S}y_{t}^{S}+0y_{q+1}^{S}=v\left( S\right) .$ Thus, $%
\tsum\limits_{i\in S}\left(
\sum_{t=1}^{q}b_{t}^{i}y_{t}^{N}+0y_{q+1}^{N}\right) \geq v\left( S\right)
,\forall S\subset N,$ and this implies that $\left( b^{i}y^{N}\right) _{i\in
N}\in C\left( v\right) .$\hfill $\blacksquare $\bigskip

\begin{corollary}
Let $\left( A,B,p,r,c\right) $ be an $LPP$ situation with $M^{\text{min}%
}=\varnothing ,$ $\left( N,\mathcal{P}(N),\left\{ V\left( \left. \bullet
\right\vert P\right) \right\} _{P\in \mathcal{P}(N)}\right) $\ the
corresponding partition function form game and $\left( N,v\right) $ the
related game in characteristic function form. Then, $C\left( V\right) =%
\overline{C}\left( V\right) =C\left( v\right) .$
\end{corollary}

Linear production $(LP)$ situations and corresponding cooperative games were
introduced in Owen (1975), where it is shown that these games have a
nonempty core by constructing a core-element via a related dual linear
program. Gellekom et al. (2000) named the set of all the core-elements that
can be found in the same way as performed by Owen, the `Owen set'. Following
this idea, we can introduce the Owen set of an $LPP$ situation $\left(
A,B,p,r,c\right) ,$ $Owen\left( A,B,p,r,c\right) ,$ as the set whose
elements can be obtained through an optimal solution of (\ref{dual N})
associated with the optimal solution of (\ref{eq:LP2}) $\left(
x^{N};d_{N}\right) $ such that $d_{N}\leq r.$ To sum up, when $M^{\text{min}%
}=\varnothing $ one way to obtain a stable distribution of the total profit
is to use an element of the so-called Owen set of the $LPP$ situation $%
\left( A,B,p,r,c\right) .$ We should mention that this is similar to the
classical results in the $LP$ situations. However, it does not always work
in the same way. Although the games such as those in Proposition \ref{N>r}
are characteristic function form games, they can have an empty core, as the
following example shows.

\begin{example}
\label{ex c vacio}Let $(A,B,r,p,c)$ be an $LPP$ situation, with three
producers, $N=\left\{ 1,2,3\right\} ,$ who produce three products from three
resources and a common-pool resource, where 
\begin{equation*}
A=\left[ 
\begin{array}{lll}
3 & 6 & 6 \\ 
6 & 6 & 6 \\ 
5 & 10 & 6 \\ 
2 & 4 & 4%
\end{array}%
\right] ,B=\left[ 
\begin{array}{ccc}
15 & 6 & 9 \\ 
4 & 18 & 9 \\ 
16 & 19 & 2%
\end{array}%
\right] ,p=\left( 
\begin{array}{c}
10 \\ 
9 \\ 
9%
\end{array}%
\right) ,c=2,r=10.
\end{equation*}%
The demands are 
\begin{equation*}
\begin{array}{l}
d_{\left\{ 1\right\} }=\frac{4}{3},d_{\left\{ 2\right\} }=4,d_{\left\{
3\right\} }=\frac{4}{5}, \\ 
d_{\left\{ 12\right\} }=\frac{22}{3},d_{\left\{ 13\right\} }=\frac{17}{3}%
,d_{\left\{ 23\right\} }=\frac{42}{5},d_{N}=\frac{31}{3},%
\end{array}%
\end{equation*}%
and $\min \left\{ \frac{31}{3},10\right\} =10.$\ Since $M^{\text{min}%
}=\left\{ N\right\} $ the corresponding TU-game $\left( N,v\right) $
associated with this situation is given by 
\begin{equation*}
\begin{array}{l}
v\left( \left\{ 1\right\} \right) =4,v\left( \left\{ 2\right\} \right)
=12,v\left( \left\{ 3\right\} \right) =\frac{12}{5}, \\ 
v\left( \left\{ 12\right\} \right) =22,v\left( \left\{ 13\right\} \right)
=13,v\left( \left\{ 23\right\} \right) =\frac{126}{5},v\left( N\right) =30,%
\end{array}%
\end{equation*}%
and $C\left( v\right) =\varnothing .$ Obviously, $C\left( V\right) =%
\overline{C}\left( V\right) =\varnothing $
\end{example}

Looking at this example one can observe that the partitions $\left\{
N\backslash \left\{ i\right\} ,\left\{ i\right\} \right\} _{i\in N}$ are the
only stable in the following sense.

\begin{definition}
A partition $P\in \mathcal{P}(N)$ is said to be partitionally stable if the
following two conditions hold $\forall S\in P$
\end{definition}

\begin{enumerate}
\item[$\left( 1\right) $] $C\left( v^{S}\right) \neq \varnothing $, \textit{%
and }

\item[$\left( 2\right) $] $\nexists \left\{ T_{k}\right\} _{k=1}^{l}\in P$%
\textit{\ such that} $C\left( v^{S\cup \left\{
\tbigcup\limits_{k=1}^{l}T_{k}\right\} }\right) \neq \varnothing ,$
\end{enumerate}

\noindent \textit{where }$\left( S,v^{S}\right) $\textit{\ is the game
reduced to coalition }$S$ when partition $P$ is formed\textit{.\bigskip }

This definition of stability holds for games in partition function form
(with whatever concept of dominance) and in characteristic function form and
allows us, in some sense, to extend the concept of the core. Note that when
the core of a game is non empty, then the grand coalition is the only one
which is partitionally stable.

\begin{proposition}
Let $\left( A,B,p,r,c\right) $ be an $LPP$ situation with $\left\{ N\right\}
\in M^{\text{min}}.$ The partitions $\left\{ N\backslash \left\{ i\right\}
,\left\{ i\right\} \right\} _{i\in N}$ are the only partitionally stable or
the grand coalition if the core of the game is non empty.
\end{proposition}

\noindent \textbf{Proof.} If the core of the game is non empty the result
holds. If it is empty, $d_{N\backslash \left\{ i\right\} }+d_{\left\{
i\right\} }\leq r,d_{N}>r.$ The games associated with $\left\{ N\backslash
\left\{ i\right\} \right\} _{i\in N}$ and $\left\{ \left\{ i\right\}
\right\} _{i\in N}$ have a non empty core by Proposition \ref{dual} and any
other partition does not satisfy condition (2) in the previous
definition.\hfill $\blacksquare $\bigskip

In the more general case, when the common-pool resource could be a
constraint for the production process for some partition$,$ each coalition
of producers will obtain an amount of common-pool resource from either
through a collaborative procedure or through a competitive mechanism$.$ With
respect to the information that they have on this process, we consider that
they do not know the way in which $r$ will be shared. The next section
tackles this situation.

\section{The common-pool resource game}

In this section we assume that producers do not know how the common-pool
resource is to be assigned. Therefore, they do not know the partition
function form game $V,$ so they face a problem under uncertainty. Thus, they
can examine the problem of the amount they will receive from different
points of view. We consider that what a coalition of producers $S$ expects
to receive from the common-pool resource can be described by the common-pool
resource game $\left( N,R\right) .$ These games are cooperative $TU-$games
in characteristic function form and can be defined following different
approaches. There are two extreme cases, depending on which point of view is
used to deal with the situation, the optimistic and the pessimistic
common-pool resource games, that will be addressed in this section. Hence, $%
R\left( S\right) $ is what coalition $S$ thinks it can guarantee from the
common-pool resource working on their own. It can be any value between those
obtained from the optimistic and pessimistic points of view.

Once a coalition of producers $S$ received its share of the common-pool
resource, $R\left( S\right) ,$ using this amount as $z_{S}(P)$ in (\ref%
{eqLP1}), for all $S\subseteq N,$ the $LPP$\ game $\left( N,v^{R}\right) $
is obtained, where $v^{R}\left( S\right) =value\left( S;R\left( S\right)
\right) .$ In this way, the game associated with the $LPP$ situation $%
(A,B,r,p,c)$ obtained from the common-pool resource game $\left( N,R\right) $
reduces to a characteristic function form game, $\left( N,v^{R}\right) ,$
since it does not depend on what the others may do.

The following theorem states a sufficient condition for the $LPP$ game $%
\left( N,v^{R}\right) $ to have a non empty core when $d_{N}>r,$\ no matter
from what point of view the common-pool resource game $\left( N,R\right) $
is defined. The case $d_{N}\leq r$ is studied in subsection \ref%
{secoptimistic}.

\begin{theorem}
\label{prop r no vacio}Let $\left( A,B,p,r,c\right) $ be an $LPP$ situation,
let $\left( N,R\right) $ be the common-pool resource game associated with
it, and $\left( N,v^{R}\right) $ the corresponding $LPP$\ game. When $%
d_{N}>r $\ if $C\left( R\right) \neq \varnothing ,$ then $C\left(
v^{R}\right) \neq \varnothing .$
\end{theorem}

\noindent \textbf{Proof. }Since $C\left( R\right) \neq \varnothing ,$ there
is $u\in \mathbb{R}^{N}$ such that $u(S)=\sum_{i\in S}u_{i}\geq R\left(
S\right) ,$ for all $S$, and $u(N)=r.$ Let $y^{\ast }$ be an optimal
solution of the dual problem of (\ref{eq:LP3}). From duality theory, we know
that $\sum_{t=1}^{q}b_{t}^{N}y_{t}^{\ast }+ry_{q+1}^{\ast }-cr=v^{R}(N).$ On
the other hand, $\forall S\subseteq N$%
\begin{equation*}
\sum_{t=1}^{q}b_{t}^{S}y_{t}^{\ast }+u\left( S\right) y_{q+1}^{\ast
}-cu\left( S\right) \geq \sum_{t=1}^{q}b_{t}^{S}y_{t}^{\ast }+R\left(
S\right) (y_{q+1}^{\ast }-c)\geq v^{R}(S),
\end{equation*}%
where the last inequality holds because $y^{\ast }$ is feasible for the dual
problem of coalition $S$ and $y_{q+1}^{\ast }>c$ since $d_{N}>r$. Thus, $%
\left( b^{i}y^{\ast }+u_{i}\left( y_{q+1}^{\ast }-c\right) \right) _{i\in
N}\in C\left( v^{R}\right) \neq \varnothing .$\hfill $\blacksquare $\bigskip

However, the opposite is not true in general as the next example shows.

\begin{example}
\label{excv+=0}Let $(A,B,r,p,c)$ be an $LPP$ situation, with two producers, $%
N=\left\{ 1,2\right\} ,$ who produce three products from two resources and a
common-pool resource, where 
\begin{equation*}
A=\left[ 
\begin{array}{lll}
1 & 0 & 1 \\ 
0 & 1 & 1 \\ 
2 & 2 & 1%
\end{array}%
\right] ,B=\left[ 
\begin{array}{ll}
4 & 1 \\ 
1 & 4%
\end{array}%
\right] ,p=\left( 
\begin{array}{l}
4 \\ 
4 \\ 
8%
\end{array}%
\right) ,c=1,r=4.
\end{equation*}%
Consider the common-pool resource game, $\left( N,R\right) $, such that $%
R\left( \left\{ 1\right\} \right) =R\left( \left\{ 2\right\} \right)
=R\left( \left\{ 12\right\} \right) =4.$ In this case, $C\left( R\right)
=\varnothing ,$\ $d_{N}=5$\ and $v^{R}\left( \left\{ 1\right\} \right)
=v^{R}\left( \left\{ 2\right\} \right) =10,v^{R}\left( N\right) =28$, thus $%
C\left( v^{R}\right) \neq \varnothing .$
\end{example}

\subsection{The optimistic approach\label{secoptimistic}}

From an optimistic point of view, a coalition of producers $S$ will obtain
its demand. The related common-pool resource game\textit{\ }$\left(
N,R^{opt}\right) $ is such that $R^{opt}\left( S\right) =\min \left\{
d_{S},r\right\} .$

Using the amount $R^{opt}\left( S\right) $ in (\ref{eqLP1}), for all $%
S\subseteq N,$ the optimistic\textit{\ }$LPP$\textit{\ game} $\left(
N,v^{opt}\right) $ is derived.

The core of this class of games can be non empty, as Example \ref{excv+=0}\
illustrates. But $C\left( v^{opt}\right) =\varnothing $ on many occasions,
as the next example shows.

\begin{example}
\label{ex cv+ vacio}Let $(A,B,r,p,c)$ be an $LPP$ situation, with three
producers, $N=\left\{ 1,2,3\right\} ,$ who produce two products from two
resources and a common-pool resource, where 
\begin{equation*}
A=\left[ 
\begin{array}{ll}
2 & 3 \\ 
3 & 2 \\ 
1 & 1%
\end{array}%
\right] ,B=\left[ 
\begin{array}{lll}
40 & 60 & 80 \\ 
60 & 40 & 50%
\end{array}%
\right] ,p=\left( 
\begin{array}{l}
50 \\ 
60%
\end{array}%
\right) ,c=14,r=50
\end{equation*}%
and $\left( N,v^{opt}\right) $ the related optimistic\textit{\ }$LPP$\textit{%
\ game}. In this case, the core of the optimistic\textit{\ }$LPP$\textit{\
game }will be all the points in $\mathbb{R}^{3}$ such that 
\begin{equation*}
\begin{array}{l}
x_{1}\geq 720,x_{2}\geq 920,x_{3}\geq 1150, \\ 
x_{1}+x_{2}\geq 1640,x_{1}+x_{3}\geq 1936,x_{2}+x_{3}\geq
2070,x_{1}+x_{2}+x_{3}=2300,%
\end{array}%
\end{equation*}%
\ but it can be seen that there is no point satisfying all the above
inequalities, then $C\left( v^{opt}\right) =\varnothing .$ Taking into
account that 
\begin{equation*}
\begin{array}{l}
d_{\left\{ 1\right\} }=20,d_{\left\{ 2\right\} }=20,d_{\left\{ 3\right\}
}=25, \\ 
d_{\left\{ 12\right\} }=40,d_{\left\{ 13\right\} }=46,d_{\left\{ 23\right\}
}=45,d_{N}=66\text{\ and }\min \left\{ 66,50\right\} =50,%
\end{array}%
\end{equation*}%
\ it is easy to check that $C\left( R^{opt}\right) =\varnothing .$
\end{example}

We will study a situation in which $C\left( v^{opt}\right) $ is non empty,
but first we need some results. The following lemma tells us that when the
common-pool resource is sufficient for the grand coalition, then the value
of the grand coalition is an upper bound for the sum of the optimistic
values in every partition.

\begin{lemma}
Let $\left( A,B,p,r,c\right) $ be an $LPP$ situation and $\left(
N,v^{opt}\right) $ the $LPP$\ optimistic game associated with it. If $%
d_{N}\leq r,$ then%
\begin{equation*}
\tsum\limits_{S\in P}value\left( S,d_{S}\right) \leq value\left(
N,d_{N}\right) =v^{opt}\left( N\right) ,\forall P\in \mathcal{P}(N).
\end{equation*}
\end{lemma}

\noindent \textbf{Proof. }We consider the linear program $\left( \ref{eq:LP2}%
\right) $ for the grand coalition $N:$

\begin{equation}
\begin{array}{ll}
\max & \sum_{j=1}^{g}p_{j}x_{j}-cz \\ 
\text{s.t:} & Ax\leq \left( 
\begin{array}{c}
b^{N} \\ 
z%
\end{array}%
\right) \\ 
& x\geq \mathbf{0}_{g},z\geq 0.%
\end{array}
\label{eq LP4}
\end{equation}%
Given $P\in \mathcal{P}(N)$, for every $S\in P$ there is an optimal solution 
$\left( x^{S};d_{S}\right) $ for the linear program (\ref{eq:LP2}). Thus, $%
\left( \left( \tsum\limits_{S\in P}x_{j}^{S}\right)
_{j=1}^{g};\tsum\limits_{S\in P}d_{S}\right) $ is a feasible solution for
the linear program (\ref{eq LP4}), therefore we have 
\begin{equation}
\sum\limits_{j=1}^{g}p_{j}\left( \sum\limits_{S\in P}x_{j}^{S}\right)
-c\left( \tsum\limits_{S\in P}d_{S}\right) \leq
\sum\limits_{j=1}^{g}p_{j}x_{j}^{N}-cd_{N},  \label{**}
\end{equation}%
where $\left( x^{N};d_{N}\right) $ is an optimal solution for (\ref{eq LP4})
such that $d_{N}\leq r.$ If we rewrite $\left( \ref{**}\right) $, we obtain $%
\sum\limits_{S\in P}\sum\limits_{j=1}^{g}\left( p_{j}x_{j}^{S}-cd_{S}\right)
\leq \sum\limits_{j=1}^{g}p_{j}x_{j}^{N}-cd_{N}$ and $\sum\limits_{S\in
P}value\left( S,d_{S}\right) \leq v^{opt}\left( N\right) .$\hfill $%
\blacksquare $\bigskip

The following lemma, which is given without a proof because it is easy to
derivate, give us two linear programs that, although they have different
optimal solution sets, also have the same optimal values, i.e. they are
optimally equivalents. Note that an optimal solution of the second one is
the optimal demand of the common-pool resource for each coalition $S,$ $%
d_{S}.$ We should highlight that they only differ in a redundant constrain, $%
z\leq d_{S}$, however, this is the key with which to prove the next theorem.

\begin{lemma}
Let $\left( A,B,p,r,c\right) $ be an $LPP$ situation. The following linear
programs are optimally equivalents, for all $S$,%
\begin{equation}
\begin{array}{ll}
\max & \sum_{j=1}^{g}p_{j}x_{j}-cz \\ 
\text{s.t:} & Ax\leq \left( 
\begin{array}{c}
b^{S} \\ 
z%
\end{array}%
\right) \\ 
& z\leq d_{S} \\ 
& x\geq \mathbf{0}_{g},z\geq 0.%
\end{array}
\label{eqLPi}
\end{equation}%
\begin{equation}
\begin{array}{ll}
\max & \sum_{j=1}^{g}p_{j}x_{j}-cz \\ 
\text{s.t:} & Ax\leq \left( 
\begin{array}{c}
b^{S} \\ 
z%
\end{array}%
\right) \\ 
& x\geq \mathbf{0}_{g},z\geq 0.%
\end{array}%
\end{equation}
\end{lemma}

The previous results provide us with two different, but equivalent, ways in
which to tackle the linear programs. In the proof of the following theorem,
we use one or the other depending on which will be more helpful. The next
result tells us that cooperation eliminates the conflict. Because they could
all go together to request the amount of the common-pool resource no matter
which mechanism the manager uses to distribute it, regardless of what
happens with the rest of partitions.

\begin{theorem}
\label{thdn<r}Let $\left( A,B,p,r,c\right) $ be an $LPP$ situation and $%
\left( N,v^{opt}\right) $ the $LPP$\ optimistic game associated with it. If $%
d_{N}\leq r,$ then $C\left( v^{opt}\right) \neq \varnothing .$
\end{theorem}

\noindent \textbf{Proof. }Consider the linear program $\left( \ref{eqLPi}%
\right) $ for the grand coalition, 
\begin{equation}
\begin{array}{ll}
\max & \sum_{j=1}^{g}p_{j}x_{j}-cz \\ 
\text{s.t:} & Ax\leq \left( 
\begin{array}{c}
b^{N} \\ 
z%
\end{array}%
\right) \\ 
& z\leq d_{N} \\ 
& x\geq \mathbf{0}_{g},z\geq 0.%
\end{array}
\label{eqlpiN}
\end{equation}%
\ its dual is given by%
\begin{equation}
\begin{array}{ll}
\min & \sum_{t=1}^{q}b_{t}^{N}y_{t}+0y_{q+1}+d_{N}y_{q+2} \\ 
\text{s.t:} & A^{t}y\geq p \\ 
& y_{q+1}-y_{q+2}\leq c \\ 
& y\geq \mathbf{0}_{q+2}.%
\end{array}
\label{eqLPdi}
\end{equation}%
Let $\left( x^{N};d_{N}\right) $ and $\left( y_{q}^{N},y_{q+1}^{N},0\right) $
be the primal and dual optimal solutions for $\left( \ref{eqlpiN}\right) $
and $\left( \ref{eqLPdi}\right) $, respectively with $d_{N}\leq r$ and $%
y_{q+2}^{N}=0\footnote{%
This is true because the common-pool resource is not scarce in this case.}.$
It is easy to check that $\left( y_{q}^{N},y_{q+1}^{N},0\right) $ is a
feasible solution for the dual problem of $\left( \ref{eqLPi}\right) $ for
every coalition $S$. If $\left( y_{q}^{S},y_{q+1}^{S},y_{q+2}^{S}\right) $
is an optimal dual solution associated with $\left( x^{S};d_{S}\right) $, it
holds that $\sum_{t=1}^{q}b_{t}^{S}y_{t}^{N}+0y_{q+1}^{N}+d_{S}y_{q+2}^{N}%
\geq
\sum_{t=1}^{q}b_{t}^{S}y_{t}^{S}+0y_{q+1}^{S}+d_{S}y_{q+2}^{S}=value(S,d_{S})=v^{opt}\left( S\right) . 
$ Therefore, $\sum\limits_{i\in S}\left(
\sum_{t=1}^{q}b_{t}^{i}y_{t}^{N}\right) \geq v^{opt}\left( S\right) ,$ $%
\forall S\subseteq N,$ and this implies that $\left( b^{i}y^{N}\right)
_{i\in N}\in C\left( v^{opt}\right) .$\hfill $\blacksquare $\bigskip

The following result is given without a proof because it is straightforward,
since $v^{opt}\geq v^{+}\geq v^{-}.$

\begin{corollary}
Let $\left( A,B,p,r,c\right) $ be an $LPP$ situation, $\left( N,\mathcal{P}%
(N),\left\{ V\left( \left. \bullet \right\vert P\right) \right\} _{P\in 
\mathcal{P}(N)}\right) $\ the corresponding partition function form game and 
$\left( N,v^{+}\right) ,\left( N,v^{+}\right) $ the related games in
characteristic function form. If $d_{N}\leq r,$ then 
\begin{equation*}
\begin{array}{l}
1.\text{ }C\left( v^{+}\right) \neq \varnothing \text{ and }\overline{C}%
\left( V\right) \neq \varnothing  \\ 
2.\text{ }C\left( v^{-}\right) \neq \varnothing \text{ and }C\left( V\right)
\neq \varnothing .%
\end{array}%
\end{equation*}
\end{corollary}

\begin{remark}
Note that this theorem holds for all $LPP$\ games $\left( N,v^{R}\right) ,$
obtained from any common-pool resource game $\left( N,R\right) $\ associated
with an $LPP$ situation, since $v^{opt}\left( S\right) \geq v^{R}\left(
S\right) ,$ for every coalition $S$.
\end{remark}

This result is important for several reasons. Firstly, we have found a case
in which the core of the optimistic game is non empty. Secondly, the Owen
set is very easy to obtain. Thirdly, it shows that cooperation among all
agents is important when it makes the common-pool resource not to be scarce.
Finally, in this case no matter how $\left( N,R\right) $ or $z_{S}\left(
P\right) $ are.

At a first glance, it seems that an easy condition to assure that the core
is empty, when $d_{N}>r,$ could be $\exists P\in \mathcal{P}(N)$ such that $%
\tsum\limits_{S\in P}d_{S}>r,$ however, is not true as the next example
shows.

\begin{example}
Let $(A,B,r,p,c)$ be the $LPP$ situation described in Example \ref{excv+=0}.
In this case, $d_{\left\{ 1\right\} }=d_{\left\{ 2\right\} }=7,d_{N}=5,$ $%
v^{opt}\left( 1\right) =v^{opt}\left( 2\right) =10$ and $v^{opt}\left(
N\right) =28.$ Thus, there is a partition $P=\left\{ \left\{ 1\right\}
,\left\{ 2\right\} \right\} $ where $d_{\left\{ 1\right\} }+d_{\left\{
2\right\} }>4$ and the core is non empty. Therefore, the aforementioned
condition does not guarantee that the core is empty.
\end{example}

When $d_{N}>r$ and, $\forall P\in \mathcal{P}(N),$ $\tsum\limits_{S\in
P}d_{S}<r$ the core of the optimistic game can be empty as Example \ref{ex c
vacio} shows. But it does not hold in general as the following example
illustrates.

\begin{example}
\label{exdn>rdi<r}Let $(A,B,r,p,c)$ be an $LPP$ situation, with three
producers, $N=\left\{ 1,2,3\right\} ,$ who produce three products from three
resources and a common-pool resource, where 
\begin{equation*}
A=\left[ 
\begin{array}{lll}
10 & 8 & 7 \\ 
7 & 10 & 5 \\ 
3 & 6 & 7 \\ 
5 & 2 & 4%
\end{array}%
\right] ,B=\left[ 
\begin{array}{lll}
9 & 6 & 8 \\ 
5 & 18 & 6 \\ 
17 & 13 & 3%
\end{array}%
\right] ,p=\left( 
\begin{array}{l}
8 \\ 
9 \\ 
5%
\end{array}%
\right) ,c=1,r=5.
\end{equation*}%
The corresponding optimistic game is 
\begin{equation*}
\begin{array}{ccc}
v^{opt}\left( 1\right) =1.5, & v^{opt}\left( 2\right) =5.25, & v^{opt}\left(
3\right) =3.5, \\ 
v^{opt}\left( 12\right) =13.125, & v^{opt}\left( 13\right) =7.7, & 
v^{opt}\left( 23\right) =12.25,%
\end{array}%
\end{equation*}%
and $v^{opt}\left( N\right) =17.5.$ The demands are%
\begin{equation*}
\begin{array}{ccc}
d_{\left\{ 1\right\} }=1, & d_{\left\{ 2\right\} }=1.5, & d_{\left\{
3\right\} }=1, \\ 
d_{\left\{ 12\right\} }=3.75, & d_{\left\{ 13\right\} }=2.2, & d_{\left\{
23\right\} }=3.5,%
\end{array}%
\end{equation*}%
with $d_{N}>5$\ and the core is non empty.
\end{example}

Hence, when $d_{N}>r$ we have from Theorem \ref{prop r no vacio}\ a
condition which is sufficient for the non emptiness of the core. In general,
it is not clear whether the core of the optimistic game is empty or not.

\subsection{The pessimistic approach}

From a pessimistic point of view, a coalition of producers $S$ will receive
what agents outside $S$ leave using the partition that minimizes the
remainder for $S$. This situation can be described as a common-pool resource
game $\left( N,R^{pes}\right) ,$ where $R^{pes}\left( S\right) =\min \left\{
\min\limits_{P:S\in P}\left\{ \left( r-\tsum\nolimits_{\substack{ T\in P  \\ %
T\neq S}}d_{T}\right) _{+}\right\} ,d_{S}\right\} .$

Using this amount $R^{pes}\left( S\right) $ as $z_{S}(P)$ in (\ref{eqLP1}),
for all $S\subseteq N,$ the pessimistic $LPP$\textit{\ game} $\left(
N,v^{pes}\right) $ is obtained.

When $d_{N}\leq r$ the core of this game is non empty since the core of the
optimistic game is non empty and $v^{opt}\geq v^{+}\geq v^{-}\geq v^{pes}$.
However, when $d_{N}>r$ it can be empty as Example \ref{ex c vacio} shows
and, therefore, $C\left( v^{+}\right) =\overline{C}\left( V\right) =C\left(
v^{-}\right) =C\left( V\right) =\varnothing $. We should point out that in
Example \ref{ex c vacio} the optimistic and pessimistic games coincide. The
following result states a condition for the non emptiness of the core of the
pessimistic game.

\begin{theorem}
Let $\left( A,B,p,r,c\right) $ be an $LPP$ situation and $\left(
N,v^{pes}\right) $ the \textit{pessimistic }$LPP$\textit{\ game} associated
with it. If $d_{N}>r$ and $\tsum\limits_{i\in N}d_{i}\geq r$, then $C\left(
v^{pes}\right) \neq \varnothing .$
\end{theorem}

\noindent \textbf{Proof. }Let $\left\{ d_{i}\right\} _{i\in N}$ be the
individual demands of agents in $N$. Consider the common-pool resource game $%
(N,w)$, where $w\left( S\right) =\left( r-\tsum\nolimits_{i\notin
S}d_{i}\right) _{+}$ and $w\left( N\right) =r.$

We will distinguish two cases:

\begin{enumerate}
\item[a) ] If $\tsum\limits_{i\in N}d_{i}>r$, $(N,w)$ is a standard
bankruptcy game and, therefore, it has a non empty core. Then, an $u\in 
\mathbb{R}^{N}$ such that $u\left( N\right) =r$ exists and $u\left( S\right)
\geq \left( r-\tsum\nolimits_{i\notin S}d_{i}\right) _{+}\geq
\min\limits_{P:S\in P}\left\{ \left( r-\tsum\nolimits_{\substack{ T\in P  \\ %
T\neq S}}d_{T}\right) _{+}\right\} \geq \min \left\{ \min\limits_{P:S\in
P}\left\{ \left( r-\tsum\nolimits_{\substack{ T\in P  \\ T\neq S}}%
d_{T}\right) _{+}\right\} ,d_{S}\right\} =R^{pes}\left( S\right) .$
\end{enumerate}

Therefore, $C\left( R^{pes}\right) \neq \varnothing $ and using the same
arguments as in Theorem \ref{prop r no vacio}, the result holds.

\begin{enumerate}
\item[b) ] If $\tsum\limits_{i\in N}d_{i}=r$, $d\left( S\right)
=\tsum\limits_{i\in S}d_{i}\geq \left( r-\tsum\nolimits_{i\notin
S}d_{i}\right) _{+}\geq \min\limits_{P:S\in P}\left\{ \left(
r-\tsum\nolimits _{\substack{ T\in P  \\ T\neq S}}d_{T}\right) _{+}\right\}
\geq \min \left\{ \min\limits_{P:S\in P}\left\{ \left( r-\tsum\nolimits 
_{\substack{ T\in P  \\ T\neq S}}d_{T}\right) _{+}\right\} ,d_{S}\right\}
=R^{pes}\left( S\right) .$
\end{enumerate}

Thus, $C\left( R^{pes}\right) \neq \varnothing $ and, then, from Theorem \ref%
{prop r no vacio} $C\left( v^{pes}\right) \neq \varnothing .$\hfill $%
\blacksquare $\bigskip

When $d_{N}>r$ and $\tsum\limits_{i\in N}d_{i}<r,$ the core of the
pessimistic game can be empty as in Example \ref{ex c vacio} or it can be
non empty as Example \ref{exdn>rdi<r}\ illustrates, since if the core of the
optimistic game is non empty, the core of the pessimistic one is also non
empty. Thus, in this case we have obtained similar results to those applying
the optimistic approach.

\section{Concluding remarks}

The model proposed in this paper is novel, arises from many real-life
situations and contains important changes with respect to linear production
situations that affect cooperation. In spite of a small change in the $LP$
model we obtain in a natural way, different situations where the
corresponding games are, in general, games in partition function form as
opposed to games in characteristic function form. Moreover, contrary to $LP$
games these games can have an empty core.

We should highlight the role that demand of the common-pool resource for the
grand coalition$,$ $d_{N},$ plays when tackling the problem. We have come
across an interesting case where cooperation makes the common-pool resource
not to be scarce$,$ $d_{N}\leq r,$ in which the core of the corresponding
game is non empty. When the common-pool resource is not sufficient$,$ $%
d_{N}>r,$ additional conditions are needed to assure the non emptiness of
the core. The study of this case when the partition function is known
exactly could be addressed using different approaches: from a
non-cooperative point of view, such as that used in Gutierrez et al (2015),
through an auction mechanism or with bankruptcy techniques which we would
like to give our attention to in future research. \bigskip

\textbf{Acknowledgement. }Financial support from the Government of Spain and
FEDER under projects MTM2011-23205, MTM2011-27731-C03, MTM2014-53395-C3-3-P
and MTM2014-54199-P, from Fundaci\'{o}n S\'{e}neca de la Regi\'{o}n de
Murcia through grant 19320/PI/14 and from Xunta de Galicia under the project
INCITE09-207-064-PR are gratefully acknowledged.

\end{document}